\documentstyle{article}
\title{Three alternating sign matrix identities in search of bijective proofs}
\author{David M.\ Bressoud}
\date{\today \\[5pt] {\it In honor of the master of bijective proofs, Dominique Foata.}}
\newcommand{\lmatrix}{\left( \begin{array}}
\newcommand{\rmatrix}{\end{array} \right)}
\newcommand{\ldet}{\left| \begin{array}}
\newcommand{\rdet}{\end{array} \right|}
\newcommand{\inv}{{\cal{I}}}
\newcommand{\calA}{{\cal{A}}}
\newcommand{\calT}{{\cal{T}}}
\newcommand{\calS}{{\cal{S}}}
\newcommand{\bin}[2]{\left( #1 \atop #2 \right)}
\newtheorem{proposition}{Proposition}
\newtheorem{problem}{Problem}

\begin{document}
\maketitle 

\begin{abstract}This paper highlights three known identities, each of which involves sums over
alternating sign matrices. While proofs of all three are known, the only known derivations are
as corollaries of difficult results. The simplicity and natural combinatorial interpretation of
these identities, however, suggest that there should be direct, bijective proofs.
\end{abstract}

\section{Introduction}

Alternating sign matrices (ASMs) are square matrices of 0s, 1s, and $-1$s with row and
column-sums equal to 1 and with the restriction that the non-zero entries alternate signs
across each row and down each column. An example is
\[ \left( \begin{array}{ccccc}
0 & 1 & 0 & 0 & 0 \\ 1 & -1 & 0 & 1 & 0 \\ 0 & 1 & 0 & -1 & 1 \\ 0 & 0 & 0
& 1 & 0 \\ 0 & 0 & 1 & 0 & 0 
\end{array} \right)
\]
These are rich combinatorial objects with connections to many
problems in algebraic combinatorics (see \cite{bre}, \cite{BP}, \cite{rob}). They also have
many different representations. The representation that was used in Kuperberg's proof of the
counting function for alternating sign matrices \cite{kup} and Zeilberger's proof of the
refined alternating sign matrix conjecture \cite{zeil} is the six-vertex model of statistical
mechanics. These are directed graphs in which each vertex has in-degree two and out-degree
two, and boundary conditions that the vertical arrows along the top and bottom are directed
out, horizontal arrows along the left and right are directed in, as in the following directed
graph.

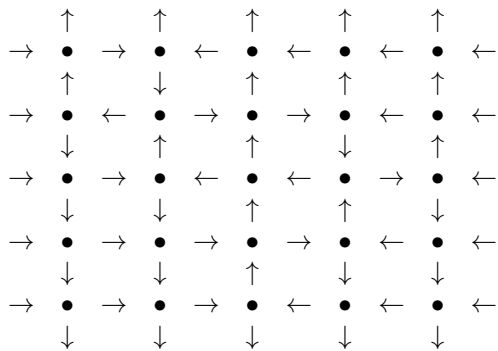
\[ \begin{array}{ccccccccccc}
 & \uparrow &  & \uparrow &  & \uparrow &  & \uparrow &  & \uparrow &
 \\
\rightarrow & \bullet & \rightarrow & \bullet & \leftarrow & \bullet &
\leftarrow & \bullet & \leftarrow  & \bullet & \leftarrow \\
 & \uparrow &  & \downarrow &  & \uparrow &  & \uparrow &  & \uparrow
&  \\
\rightarrow & \bullet & \leftarrow & \bullet & \rightarrow & \bullet &
\rightarrow & \bullet & \leftarrow  & \bullet & \leftarrow \\
 & \downarrow &  & \uparrow &  & \uparrow &  & \downarrow &  &
\uparrow & 
\\
\rightarrow & \bullet & \rightarrow & \bullet & \leftarrow & \bullet &
\leftarrow & \bullet & \rightarrow  & \bullet & \leftarrow \\
 & \downarrow &  & \downarrow &  & \uparrow &  & \uparrow &  &
\downarrow &  \\
\rightarrow & \bullet & \rightarrow & \bullet & \rightarrow & \bullet &
\rightarrow & \bullet & \leftarrow  & \bullet & \leftarrow \\
 & \downarrow &  & \downarrow &  & \uparrow &  & \downarrow &  &
\downarrow &  \\
\rightarrow & \bullet & \rightarrow & \bullet & \rightarrow & \bullet &
\leftarrow & \bullet & \leftarrow  & \bullet & \leftarrow \\
 & \downarrow &  & \downarrow &  & \downarrow &  & \downarrow &  &
\downarrow & 
\end{array} \]

To actually make this a directed graph on 25 vertices, we can identify the $i$th up-arrow along
the top row with the $i$th right arrow along the left edge, and similarly identify bottom and
right arrows. This is called a six-vertex model because there are six possible configurations
at each vertex. We shall describe a vertex as {\bf horizontal} if both in-edges are
horizontal, {\bf vertical} if both in-edges are vertical, and otherwise {\bf southwest}, {\bf
northwest}, {\bf northeast}, or {\bf southeast}, according to the direction of the sum of the
four vectors represented by the four adjacent edges. 

It should be noted that the sum of all vertical vectors is zero, as is the sum of all
horizontal vectors. It follows that there will always be an equal number of southwest and
northeast vertices, and an equal number of southeast and northwest vertices.

Our example of a six-vertex model corresponds to our example of an alternating
sign matrix. Each 1 in the ASM corresponds to a horizontal vertex, each $-1$ to a vertical
vertex, and the 0s to the other vertices. This is a bijection because once the
positions of the horizontal and vertical vertices are known, all other vertices are uniquely
determined.

The six-vertex model is not the only insightful representation, but it is very suggestive,
especially because there is also a natural connection between ASMs and complete directed
graphs or tournaments. It would be very useful to have a direct bijective connection between
ASMs and tournaments. In explaining the bijection that we seek, we shall also present two
other related identities that cry out for bijective proofs. 

\section{The $\lambda$-determinant}

The first two identities that I wish to present arise from the $\lambda$-determinant of Robbins
and Rumsey \cite{RR}. This is based on the Desnanot--Jacobi adjoint matrix theorem \cite{des},
\cite{jac} that was used by Dodgson \cite{dodg} to create his algorithm for evaluating
determinants. Given a square matrix $M$, we let $M^i_j$ denote $M$ with row $i$ and column $j$
deleted. We then have that
\begin{equation}
\det M\ =\ \frac{\det M_1^1 \cdot \det M_n^n - \det M_1^n \cdot \det M_n^1}{\det
M_{1,n}^{1,n}}. \label{eqn:1}
\end{equation}

If we define the determinant of an empty matrix ($0\times 0$) to be 1 and the determinant of
the $1\times 1$ matrix $(a)$ to be $a$, then equation~(\ref{eqn:1}) can be used as a recursive
definition of the determinant. A natural one-parameter generalization of the determinant
arises if we use the same initial conditions and replace the minus sign in the numerator of
the recursive step by $+\lambda$:
\begin{equation} {\det}_{\lambda}(M) = \frac{\displaystyle {\det}_{\lambda}\left(M_1^1\right)
{\det}_{\lambda}\left(M_n^n\right) + \lambda\,
{\det}_{\lambda}\left(M_1^n\right)
{\det}_{\lambda}\left(M_n^1\right) }{\displaystyle
{\det}_{\lambda}\left(M_{1,n}^{1,n}
\right)}. \label{eqn:2} \end{equation}
The following generalization of the Vandermonde determinant evaluation follows by induction.

\begin{proposition}
\begin{equation} {\det}_{\lambda} \left(x_i^{n-j}\right)\ =\ \prod_{1 \leq i < j \leq n} (x_i +
\lambda x_j). \label{eqn:3} \end{equation}
\end{proposition}

If we expand a few $\lambda$-determinants, an interesting pattern emerges:
\begin{eqnarray*}&& \hspace{-10mm}{\det}_{\lambda} \lmatrix{ccc} a & b & c \\
d & e & f
\\ g & h & i\rmatrix \\[10pt]
& = & aei + \lambda\,(bdi + afh) + \lambda^2\,(bfg + cdh) + \lambda^3\,ceg +\
\lambda(1+\lambda) bde^{-1}fh, \\[15pt] && \hspace{-10mm} {\det}_{\lambda} \lmatrix{cccc} a & b
& c & d \\
 e & f & g & h
\\ i & j & k & l \\ m & n & o & p \rmatrix \\[10pt]
& = & \cdots + \lambda^3(1+\lambda) bef^{-1}hkn + \lambda^3(1+\lambda)^2
cfg^{-1} hij^{-1}kn + \cdots .
\end{eqnarray*}
The monomials in roman letters that correspond to permutation matrices are each multiplied by
$\lambda$ raised to the inversion number of the permutation. The other monomials in roman
letters that appear, such as $cfg^{-1} hij^{-1}kn$, correspond to alternating sign matrices,
in this case
\[ \lmatrix{cccc} 0 & 0 & 1 & 0 \\ 0 & 1 & -1 & 1 \\ 1 & -1 & 1 & 0 \\ 0 & 1 & 0 & 0 \rmatrix.
\]
Each of these monomials is multiplied by a power of $\lambda$ and a power of $1 + \lambda$. 

Let $\calA_n$ be the set of $n \times n$ ASMs. Given $A = (a_{ij}) \in \calA_n$,
we define its {\bf inversion number},
$\inv(A)$, to be
\[ \inv(A)\ =\ \sum_{i<k,\ j>l} a_{ij} \cdot a_{kl}. \] 
We define $N(A)$ to be the number of $-1$s in $A$. The following characterization of the
$\lambda$-determinant was published by Robbins and Rumsey in 1986 \cite{RR}.

\begin{proposition}
\begin{equation} {\det}_{\lambda}(m_{ij})\ = \ \sum_{A\in \calA_n} \lambda^{\inv(A) - N(A)} (1
+
\lambda)^{N(A)} \prod_{i,j=1}^n m_{ij}^{a_{ij}}.
\end{equation}
\end{proposition}

Zeilberger \cite{zeil2} has given a bijective proof of equation~(\ref{eqn:1}). It would be
desirable to have a direct proof of Proposition~2 by finding a similar proof of
equation~(\ref{eqn:2}) when the $\lambda$-determinant is defined by the right side of
Proposition~2.

\begin{problem}
Find a direct, bijective proof of the following identity. Within each summation, the range of
indices for the alternating sign matrices $B$ and $C$ is specified by the
product term.
\begin{eqnarray*} && \hskip -15mm \sum_{(B,C)\in\calA_n\times\calA_{n-2}}
\lambda^{\inv(B)+\inv(C)-N(B)-N(C)}(1+\lambda)^{N(B)+N(C)}
\prod_{i,j=1}^nm_{ij}^{b_{ij}} \prod_{i,j=2}^{n-1} m_{ij}^{c_{ij}} \\[20pt]
& = & \sum_{(B,C)\in\calA_{n-1}\times\calA_{n-1}}
\lambda^{\inv(B)+\inv(C)-N(B)-N(C)}(1+\lambda)^{N(B)+N(C)} \\[10pt]
&& \hskip 1in \times\ \prod_{1 \leq i < n \atop 1 \leq j < n}
m_{ij}^{b_{ij}}
\prod_{1 < i \leq n \atop 1 < j \leq n}  m_{ij}^{c_{ij}}
\\[10pt]
&& +\ \lambda \sum_{(B,C)\in\calA_{n-1}\times\calA_{n-1}}
\lambda^{\inv(B)+\inv(C)-N(B)-N(C)}(1+\lambda)^{N(B)+N(C)} \\[10pt]
&& \hskip 1in \times\ 
\prod_{1 \leq i < n \atop 1 \leq j < n} m_{ij}^{b_{ij}}
\prod_{1 < i \leq n \atop 1 < j \leq n}  m_{ij}^{c_{ij}}.
\end{eqnarray*}
\end{problem}

\section{Directed Graphs}

If we combine Propositions~1 and 2, we get that
\begin{equation} \prod_{1 \leq i < j \leq n} (x_i +
\lambda x_j)\ =\ \sum_{A\in \calA_n} \lambda^{\inv(A) - N(A)} (1
+
\lambda)^{N(A)} \prod_{i,j=1}^n x_i^{(n-j)a_{ij}}
\label{eqn:5} \end{equation}
It is worth noting that analogs of this identity for other root systems have been found by
Okada \cite{oka}. 

The left side of equation~(\ref{eqn:5}) can be interpreted as a sum over the set of
tournaments on $n$ vertices, $calT_n$. Each binomial $x_i + \lambda x_j$ corresponds to the
edge between vertices
$i$ and $j$. If the edge is directed from $i$ to $j$, we choose $x_i$. If it is directed from
$j$ to $i$, we choose $\lambda x_j$. Each tournament corresponds to a monomial in which the
power of $x_i$ is $\omega(i)$, the out-degree of vertex $i$, and the power of $\lambda$ is
$U(T)$, the number of {\bf upsets} in the tournament: $j>i$ and $j \rightarrow i$:
\[ \prod_{1 \leq i < j \leq n} (x_i +
\lambda x_j)\ =\ \sum_{T\in\calT_n} \lambda^{U(T)} \prod_{i=1}^n x_i^{\omega(i)}.
\]

We shall use the six-vertex model to interpret the right side of equation~(\ref{eqn:5}). We
begin with the following observations which are explained below.

\begin{proposition} Let $A$ be an $n \times n$ ASM. In the corresponding six-vertex
model\begin{itemize}
\item the number of horizontal vertices is $n + N(A)$,
\item the number of vertical vertices is $N(A)$,
\item the number of southwest or northeast vertices is $\inv(A) - N(A)$,
\item the number of southeast or northwest vertices is $\bin{n}{2} - \inv(A)$.
\end{itemize}
\end{proposition}

The number of vertical vertices is immediate from the bijection, and there most be one more 1
than $-1$ in each row. A southwest vertex corresponds to a 0 of the ASM for which there is a 1
above it in its column (due north) with no other non-zero entries in between, and a 1 to its
left in its row (due west) with no other non-zero entries in between. 
\[ \begin{array}{ccccc}
&&&& 1 \\ &&&& 0 \\ &&&& \vdots \\ &&&& 0 \\ 1 & 0 & \cdots & 0 & 0 (\rm SW)\ {\rm or}\ -1
\end{array} \]
The inversion number is
the number of such pairs of 1's: pairs of 1's for which there are only 0s in the positions
that are both due east of the lower 1 and strictly south and west of the upper 1, and there are
only 0s in the positions that are both due south of the upper 1 and strictly north and east of
the lower 1. The entry in the unique position due east of the lower 1 and due south of the
upper 1 must be either a 0, corresponding to a southwest vertex, or a $-1$. The remaining
observations follow from the equality of the number of southwest and northeast vertices, the
equality of the number of southeast and northwest vertices, and the fact that there are
$n^2$ vertices in all.

If we let $SW(A)$, $SE(A)$, and $V(A)$ denote, respectively, the number of southwest,
southeast, and vertical vertices in $A$ and $SW_i(A)$, $SE_i(A)$, and $V_i(A)$
the number of southwest, southeast, or vertical vertices, respectively, in column $i$ of $A$,
then the right side of equation~(\ref{eqn:5}) can be written as
\[ \sum_{A\in\calA_n} \lambda^{SW(A)} (1 +
\lambda)^{V(A)} \prod_{i=1}^n x_i^{SW_i(A)+SE_i(A)+V_i(A)}
\]
Equation~(\ref{eqn:5}) is equivalent to
\begin{equation}
\sum_{T\in\calT_n} \lambda^{U(T)} \prod_{i=1}^n x_i^{\omega(i)}\ =\ \sum_{A\in\calA_n} \lambda^{SW(A)} (1 +
\lambda)^{V(A)} \prod_{i=1}^n x_i^{SW_i(A)+SE_i(A)+V_i(A)}
\label{eqn:6}
\end{equation}

This suggests a natural bijection between tournaments on $n$ vertices and six-vertex models on
$n^2$ vertices in which we have chosen a direction (left or right) at each vertical vertex.
Each vertex in the six-vertex  model that has an in-edge from the north will define an out-edge
of the tournament. Call this vertex of the six-vertex model an {\bf initiating vertex}. If an
initiating vertex is southwest, there is an out-edge to the left, and the corresponding edge
in the tournament will contribute to the upset number. If the initiating vertex is southeast,
there is an out-edge to the right, and the corresponding edge in the tournament will not
contribute to the upset number. If the initiating vertex is vertical, we have a choice of
taking either the left or right out-edge. The left choice contributes one to the upset
number of the tournament; the right choice contributes nothing.

\begin{problem}
Find a bijective proof of equation~(\ref{eqn:6}).
\end{problem}

\section{The Izergin-Korepin Determinant Evaluation}

Kuperberg's proof of the alternating sign matrix conjecture and Zeilberger's proof of the
refined conjecture rest on the following determinant evaluation of Izergin \cite{Iz},
described in Korepin, Bogoliubov, and Izergin's {\it Quantum Inverse Scattering Method\/}
\cite{KBI}.

\begin{proposition}
Given $A\in \calA_n$, let $(i,j)$ be the vertex in row $i$, column $j$ of the corresponding
six-vertex model, and let $H$, $V,$ $SE,$ $SW,$ $NE,$ $NW$ be, respectively, the sets of
horizontal, vertical, southeast, southwest, northeast, and northwest vertices. For
indeterminants
$a$, $x_1, \ldots, x_n$, and $y_1,\ldots,y_n$, we have that
\begin{eqnarray}
&& \hspace{-15mm} \det \left( \frac{1}{(x_i+y_j)(ax_i+y_j)}\right)
\frac{\prod_{i,j=1}^n (x_i+y_j)(ax_i+y_j)}{\prod_{1 \leq i < j \leq n}
(x_i-x_j)(y_i-y_j)}  \nonumber \\[5pt]
& = & \sum_{A\in \calA_n} (-1)^{N(A)} (1-a)^{2N(A)} a^{\bin{n}{2} - 
\inv(A)} \nonumber \\ && \quad \times\ 
\prod_{(i,j)\in V} x_i y_j \prod_{(i,j)\in NE\cup SW} (ax_i+y_j) \prod_{(i,j)\in NW\cup
SE} (x_i+y_j). \label{eqn:7}
\end{eqnarray}
\end{proposition}

As Lascoux has pointed out \cite{las}, the right way to understand this identity is as an
extension of Cauchy's
\begin{equation}
\det \left( \frac{1}{(x_i+y_j)}\right)\prod_{i,j=1}^n
(x_i+y_j)
\prod_{1 \leq i < j \leq n}
(x_i-x_j)^{-1} (y_i-y_j)^{-1} = 1.
\label{eqn:8} \end{equation}
This is true by inspection. The determinant times the product over $i,j$ is an alternating
polynomial in the $x_i$ and in the $y_j$. Since any alternating polynomial is divisible by the
Vandermonde product, the left side of this equality is a symmetric polynomial in the $x_i$, and
it is a symmetric polynomial in the $y_j$. The degree in $x_1$ of this polynomial is zero, and
the constant can be checked by induction.

Applying this same reasoning to the left side of equation~(\ref{eqn:7}), we see that it is a
symmetric polynomial in the $x_i$ and in the $y_j$. Its degree in $x_1$ is $n-1$. On the
right, we also have a polynomial in $x_1$ of degree $n-1$. We need only check that these two
sides agree for $n$ values of $x_1$. By induction, they agree at $x_1 = -y_1/a$. If we can show
that the right side is symmetric in the $y_j$, then the identity is proven.

Symmetry follows from Baxter's triangle-to-triangle relation which was used by Izergin to
prove that
\begin{eqnarray*} && \hskip -15mm \sum_{A\in\calA_n}\prod_{(i,j)\in H} x_i(1-a) \prod_{(i,j)\in
V} (-y_j)(1-a) \\ && \times\ 
\prod_{(i,j)\in NE\cup SW} (ax_i-y_j) \prod_{(i,j)\in NW\cup
SE} (x_i-y_j)a^{1/2}
\end{eqnarray*}
is symmetric in the $x_i$, and it is symmetric in the $y_j$.

Among the corollaries of Proposition~4, we can set $a=1$ to get Borchardt's \cite{bor}
permanent-determinant identity:
\begin{equation}
\det\left(\frac{1}{(x_i+y_j)^2}\right) \frac{\prod_{i,j=1}^n (x_i+y_j)^2}{\prod_{1\leq i < j
\leq n} (x_i-x_j)(y_i-y_j)} =  {\rm perm}\left(\frac{1}{x_i+y_j}\right) \prod_{i,j=1}^n
(x_i+y_j),
\end{equation}
where
\[ {\rm perm}(a_{ij})\ :=\ \sum_{\sigma\in \calS_n} \prod_{i=1}^n a_{i,\sigma(i)}. \]
We can set $a = \omega := e^{2\pi i/3}$, $x_j = -\omega q^{j},$ and $y_j = q^{1-j}$, evaluate
the determinant, and then take the limit as $q \to 1$ to get the number of ASMs of a given
size:
\begin{equation}
\prod_{j=0}^{n-1} \frac{(3j+1)!}{(n+j)!}\ =\ \left| \calA_n \right|.
\end{equation}

If we set $a=-1$, then the matrix for which we take the determinant is $(1/(x_i^2-y_j^2))$,
which can be evaluated using Cauchy's formula, equation~(\ref{eqn:8}). The left side of
equation~(\ref{eqn:7}) simplifies to
\[ (-1)^{n(n-1)/2} \prod_{1 \leq i < j \leq n} (x_i+x_j) (y_i+y_j). \] 
From equation~(\ref{eqn:6}), each of these Vandermonde-type products can be written as a sum
over alternating sign matrices. We let Ein$_i(A)$ be the number of vertices in row $i$ with an
in-edge from the left, Nin$_j(A)$ be the number of vertices in column $j$ with an in-edge from
below. Replacing
$y_j$ by
$-y_j$ and multiplying each side by
$x_1\cdots x_n$, the case $a=-1$ is equivalent to the identity
\begin{eqnarray}
&& \hskip -15mm \sum_{(B,C)\in \calA_n \times \calA_n} 2^{N(B)+N(C)} \prod_{i=1}^n x_i^{{\rm
Ein}_i(B)}
\prod_{j=1}^n y_j^{{\rm Nin}_j(C)}  
 \nonumber \\
& = & \sum_{A\in\calA_n}(-1)^{\inv(A)-N(A)} 4^{N(A)} \prod_{(i,j)\in H} x_i \prod_{(i,j)\in
V} y_j \prod_{(i,j)\in NE} (x_i+y_j)   \nonumber \\ && \times\
\prod_{(i,j)\in SW} (-x_i-y_j) \prod_{(i,j)\in NW} (-x_i+y_j)\prod_{(i,j)\in SE} (x_i-y_j)
\label{eqn:11} \end{eqnarray}

\begin{problem}
Find a bijective proof of equation~(\ref{eqn:11}). 
\end{problem}


\begin{thebibliography}{99}

\bibitem{bor} Borchardt. Bestimmung der symmetrischen Verbindungen ihrer erzeugenden Funktion. 
{\it Journal f\"ur die reine und angewandte Mathematik\/}. {\bf 53} (1855): 193--198.

\bibitem{bre} Bressoud, David M. {\it Proofs and Confirmations: The Story of the Alternating
Sign Matrix Conjecture\/}.  Cambridge: Cambridge University Press. 1999.

\bibitem{BP} Bressoud, David M., and James Propp. How the alternating sign matrix conjecture
was solved. {\it Notices of the AMS\/} {\bf 46} (1999): 637--645.

\bibitem{des} Desnanot, P. {\it Compl\'ement de la th\'eorie des \'equations du
premier degr\'e\/}. Private publication. Paris. 1819. Described in Thomas Muir. {\it The
Theory of Determinants in the Historical Order of Development, Vol.\ I\/}. London: Macmillan
and Co. 1906.

\bibitem{dodg} Dodgson, Charles L.  Condensation of determinants. {\it Proceedings
of the Royal Society, London\/} {\bf 15} (1866): 150--155.

\bibitem{Iz} Izergin, Anatoli G. Partition function of a six-vertex model in a
finite volume. (Russian)  {\it Dokl.\ Akad.\ Nauk SSSR\/} {\bf 297} (1987): 331--333.

\bibitem{jac} Jacobi, C.\ G.\ J.  De binis quibuslibet functionibus homogeneis secundi
ordinis per substitutiones lineares in alias binas transformandis. {\it Journal fur die Reine und
Angewandt Mathematik\/}. {\bf 12} (1833): 1--69. Reprinted in {\it C.\ G.\ J.\
Jacobi: Gesammelte Werke\/}. Vol. 3, pp.\ 191--268. Berlin: Georg Reimer, 1884.

\bibitem{KBI} Korepin, V.\ E., N.\ M.\ Bogoliubov, and A.\ G.\ Izergin. {\it
Quantum Inverse Scattering Method and Correlation Functions\/}. Cambridge:
Cambridge University Press. 1993.

\bibitem{kup} Kuperberg, Greg. Another proof of the alternating sign
matrix conjecture. {\it International Mathematics Research Notes\/}
{\bf 1996}: 139--150.

\bibitem{las} Lascoux, Alain. Square-ice enumeration. The Andrews Festschrift (Maratea, 1998). 
{\it S\'eminaire Lotharingien de Combinatoire\/} {\bf 42} (1999), Art. B42.

\bibitem{oka} Okada, Soichi. Alternating sign matrices and some deformations of
Weyl's denominator formulas. {\it Journal of Algebraic Combinatorics\/} {\bf
2} (1993): 155--176.

\bibitem{rob} Robbins, David P. The story of 1, 2, 7, 42, 429, 7436,
$\ldots$. {\it The Mathematical Intelligencer\/} {\bf 13} (1991): 12--19.

\bibitem{RR} Robbins, David P., and Howard Rumsey. Determinants and alternating
sign matrices. {\it Advances in Mathematics\/} {\bf 62} (1986): 169--184.

\bibitem{zeil} Zeilberger, Doron. Proof of the refined alternating
sign matrix conjecture. {\it New York Journal of Mathematics\/} {\bf 2} (1996):
59--68. 

\bibitem{zeil2} ---------. Dodgson's determinant-evaluation rule proved by two-timing men and
women. The Wilf Festschrift (Philadelphia, PA, 1996). 
{\it Electron. J. Combin.\/} {\bf 4} (1997), no. 2, Research Paper 22.

\end{thebibliography}
\end{document}